\magnification = 1200
\hfuzz=10pt
\hsize=4.8in
\vsize=7.3in
\baselineskip=18pt
\hoffset=0.35in
\voffset=0.1in
\parindent=3pt
\def\vsni{\vskip 0.2cm}

\def\R{I\!\!R}

\def\C{I\!\!\!\!C}
\def\N{I\!\!N}

\def\O{\Omega}

\def\o{\omega}

\def\zg{\N_{0}}
\def\limsup{\mathop{\overline{\rm lim}}}

\def\s{\sigma}

\centerline{\bf ON THE RATE OF CONVERGENCE TO EQUILIBRIUM}
\centerline{\bf FOR COUNTABLE ERGODIC MARKOV CHAINS}
\vglue 0.2cm
\vglue 0.2cm\centerline{\rm Stefano Isola}
\vglue 0.4cm
\centerline{\it Dipartimento di Matematica e Informatica 
dell'Universit\`a degli Studi di Camerino}
\centerline{\it and INFM, via Madonna delle Carceri, I-62032 Camerino, Italy.}
\centerline{\it e-mail: stefano.isola@unicam.it}
\vskip 3cm
{\bf Abstract.} Using elementary methods, we prove 
that for a countable Markov
chain $P$ of ergodic degree $d > 0$ the rate of convergence 
towards the stationary distribution is subgeometric of order $n^{-d}$,
provided the initial distribution satisfies certain
conditions of asymptotic decay. An example, modelling
a renewal process and providing
a markovian approximation scheme in dynamical system theory,
is worked out
in detail, illustrating the relationships between convergence
behaviour, analytic properties of the generating functions
associated to transition probabilities
and spectral properties of the Markov operator 
$P$ on the Banach space
$\ell_1$. Explicit conditions allowing to obtain 
the actual asymptotics for the rate of 
convergence are also discussed.
\vfill \eject
\noindent
{\bf 0. INTRODUCTION.}
\vsni
\noindent
Let $S$ be a countable set and $P:S\times S \to [0,1]$ 
be a transition probability matrix.
With no loss we may set $S=\N$.
We shall assume that 
$P$ governs an irreducible, recurrent and aperiodic Markov chain 
$X=(x_n)_0^{\infty}$ with state space 
$S$. To be more precise, we set $\zg :=\N \cup \{0\}$ and
let $\O$ denote the subset of $S^{\zg}$
given by all sequences $\o=(\o_i)_{i\in \zg}$ which satisfy for any
integer $i$: $p_{\o_i\o_{i+1}}\equiv P(\o_i,\o_{i+1})>0$.
For any $n\in \zg$
we let $x_n$ be the projection on the 
$n^{\rm th}$ coordinate, i.e. $x_n(\o)=\o_n$.
Let moreover
${\bf P}_{\nu}$ be the 
probability measure with initial distribution ${\nu}$ 
(that of $x_0$) on $\O$, i.e. 
$$
{\bf P}_{\nu}\{x_n(\o)=j\}=\sum \nu_i p^{n}_{ij}=
(\nu P^n)_j\eqno(0.1)
$$
where $p^{n}_{ij}\equiv P^{(n)}(i,j)$.
Our sample space will be $\O$ equiped with the
restriction of the product $\s$-field and with probability
measure ${\bf P}_{\nu}$ for some initial distribution 
$\nu$. We shall denote by
${\bf E}_{\nu}$ the expectation w.r.t. ${\bf P}_{\nu}$.
In particular, if $\nu =\delta_i$, where $i$ is some
reference state chosen from the outset, we have
$$
{\bf P}_{i}\{x_n(\o)=j\}={\bf P}\{x_n(\o)=j|x_0(\o)=i\}=p^{n}_{ij}.\eqno(0.2)
$$
\noindent
Let $m_{i}$ be the ${\bf P}_i$-expectation 
of $\min\{n\in \zg,\, x_n(\o)=i\}$, the time of the first visit at $i$.
It is well known (see e.g. [Chu]) that if $P$ is irreducible and aperiodic then
$$
\lim_{n\to \infty} p^{n}_{ij} = {1\over m_i},\eqno(0.3)
$$
where the r.h.s. is taken to be zero in the transient and null recurrent cases when 
$m_i=\infty$. If instead $m_i$ is
finite for some 
(and hence for all) $i\in S$ then $P$ is called ergodic, or positive recurrent, 
and there is
a (unique)
probability distribution $\pi$  
on $S$ given by
${\bf \pi} = (\pi_i)_1^{\infty}=(1/m_{i})_1^{\infty}$ which is a solution 
to ${\bf \pi}={\bf \pi} P$ and thus defines a 
stationary distribution. This paper is devoted to the study of the rate of convergence
in (0.3) for ergodic chains and more generally to the rate convergence of a given initial 
distribution $\nu$ to the stationary distribution $\pi$. It is divided into two main parts. 
In the first part (Sections 1 and 2) general convergence results are stated and proved,
which relate the rate of convergence 
to a parameter $d$ called the {\sl ergodic degree}.
Roughly speaking, the ergodic degree controls in a 
continuous fashion the number of finite moments
possessed by the time of the first visit at a given state $i\in S$ (see Defintion 1).
The fact that the speed of convergence for countable state Markov chains is connected
to the number of moments of first passage times has been put forward by several 
works starting with Feller [Fe1].
In particular, using the technique of coupling, Pitman [Pi] proved that 
if first passage times have finite $r$-th moment, with $r$ a given positive integer
with $r\geq 2$, then the rate of convergence
in (0.3)
is $o(n^{-(r-1)})$. 
For other results of the same nature we refer to [Pop] and [TT].
In Theorem 1 stated below an improvement of the above results is achieved 
in that
for any real positive value of the ergodic degree $d$, 
which is assumed to be finite,
it is possible to prove subgeometric convergence to
equilibrium of order $n^{-d}$. This amounts to obtaining
subgeometric {\sl lower} bounds as well, which are here proved using 
elementary generating functions techniques.
\noindent
The relevance of obtaining sharp bounds is further 
discussed in Section 3, where 
an example modelling
a renewal process
is worked out
in detail using a different (although similar in spirit) 
method which makes use of matrix-valued 
analytic functions and allows to further sharpen
the general results of Section 1 under suitable conditions.
The main motivation is that of illustrating 
the relationships between convergence
behaviour, analytic properties of the generating functions
associated to transition probabilities
and spectral properties of the Markov operator 
$P$ on the Banach space
$\ell_1$. A second motivation is discussed in the Appendix
and comes from the fact that this example
provides
a markovian approximation scheme in dynamical system theory,
where the question of obtaining sharp subgeometric
bounds for the decay of correlations appears to be particularly
relevant (see [Is1] and [Sa]).
\noindent

\vfill\eject
\noindent
{\bf 1. ERGODIC DEGREE AND GENERAL CONVERGENCE RESULTS.}
\vsni
\noindent
In the sequel we identify sequences 
${\bf \nu} =(\nu_i)_1^{\infty} \in \ell_1(S)$, 
the corresponding row vectors 
${\bf \nu}=(\nu_1,\nu_2,\dots)$, and finite signed measures on $S$, and
define
$$
\Vert {\bf \nu} \Vert = \sum_{i=1}^{\infty}|\nu_i|.
$$
A signed measure $\nu$ satisfying $\sum \nu_l =1$ will be called
a signed distribution.
Similarly, we shall identify sequences 
${\bf u}=(u_i)_1^{\infty} \in \ell_{\infty}(S)$, the corresponding
column vectors ${\bf u}=(u_1,u_2,\dots)^t$, and bounded
functions on $S$.

\noindent
We introduce the classical {\it taboo} quantities:
$$\eqalign{
f^n_{ij}&={\bf P}\{x_l(\o) \neq j, 0<l<n, x_n(\o)=j\,
|\,x_0(\o)=i \},\cr
{ }_kp^n_{ij}&={\bf P}\{x_l(\o) \neq k, 0<l<n, x_n(\o)=j\,
|\,x_0(\o)=i \},\cr
f^*_{ij}&=\sum_{n=1}^{\infty}f^n_{ij},\qquad
{ }_kp^*_{ij}=\sum_{n=1}^{\infty}{ }_kp^n_{ij}. \cr}
$$
Clearly ${ }_jp^n_{ij}=f^n_{ij}$.
Since we have a unique 
recurrent class, $f^*_{ij}=1$
for all $i,j\in S$. Moreover, for an ergodic chain 
we have ([Chu], Chap. I.9, Thm. 5)
$$
\lim_{n\to \infty}{\sum_{k=0}^np_{ij}^n\over\sum_{k=0}^np_{ii}^n}
={\pi_j\over \pi_i}={ }_ip^*_{ij}.\eqno(1.1)
$$
The last quantity can also be viewed as the ${\bf P}_i$-mean 
number of visits to the state $j$ before return to $i$.
The relation between the $f_{ij}$'s and the transition probabilities $p_{ij}$
is given by ([Chu], Chap. I.5, Thm. 2)
$$
p_{ij}^0=\delta_{ij},\quad f^0_{ij}=0
$$
and 
$$
p_{ij}^n=\sum_{k=1}^{n}f^k_{ij}p^{n-k}_{jj}.\eqno(1.2)
$$
We also have
$$
\pi_i={1\over \sum_{n=1}^{\infty}n\,f^n_{ii}}.\eqno(1.3)
$$
For $k\in \N$, $i\in S$, let $t_k^{(i)}$ be the time of the
$k$-th entrance into state $i$, and let
$$
r_k^{(i)}(\o)=t_{k+1}^{(i)}-t_k^{(i)},\;\; k\geq 0\eqno(1.4) 
$$ 
be the sequence of times between returns
(set $t_0^{(i)}=0$). Clearly we have $r_0^{(i)}\geq 0$
and $r_k^{(i)}>0$ for $k\geq 1$.
Moreover,
the state $i$
being recurrent, $r_1^{(i)},r_2^{(i)},\dots$ 
are i.i.d. random variables under the probability ${\bf P}_i$. 
Their common distribution is given by
$$
{\bf P}_i\{r_k^{(i)}(\o)=n\}=f^n_{ii}\qquad k\geq 1.\eqno(1.5)
$$
On the other hand, having fixed an initial distribution
$\nu$ and a reference state $i$, the random variable $r_0^{(i)}$
(the delay
in the embedded renewal process) is distributed according to 
${\bf P}_{\nu}$. More specifically,
$$
{\bf P}_{\nu}\{r_0^{(i)}=n\}=\nu_i \, \delta_{n0}+
\sum_{l\neq i}\nu_lf^n_{li}.
$$
For $\gamma \geq 0$, $i,j\in S$ (and $k\geq 1$), we set
$$
M_{ij}^{(\gamma)}:={\bf E}_i(|r_k^{(j)}|^{\gamma})=
\sum_{n=1}^{\infty} n^{\gamma}f^n_{ij}.\eqno(1.6)
$$ 
Notice that $m_i\equiv M_{ii}^{(1)}$.
Given a signed distribution
$\nu$ on $S$, we also set,
$$
{ M}_{\nu i}^{(\gamma)}:={\bf E}_{\nu}(|r_0^{(i)}|^{\gamma})=
\sum_{l\neq i}\nu_l\sum_{n=1}^{\infty} n^{\gamma}\, f^n_{li}
=\sum_{l\neq i}\nu_l\, M_{li}^{(\gamma)} .
\eqno(1.7)
$$
The next result extends ([KSK], Thm. 9.65) to arbitrary (i.e. not necessarily integer) $\gamma$-values.
\vsni
\noindent
{\bf Lemma 1.} {\sl If $\gamma \geq 0$, then
$M_{ii}^{(\gamma+1)}<\infty$ if and only if
${ M}_{\pi i}^{(\gamma)}<\infty$.}
\vsni
\noindent
{\sl Proof.} 
Using the last identity in (1.1) and the decomposition
${ }_ip^{n+m}_{ii}=\sum_{l\neq i} { }_ip^m_{il} { }_ip^n_{li}$
we get
$$
\eqalign{
{ M}_{\pi i}^{(\gamma)}&=
\sum_{n=1}^{\infty}n^{\gamma}\sum_{l\neq i}\pi_l f^n_{li}
=\pi_i\sum_{n=1}^{\infty} n^{\gamma}\sum_{l\neq i} { }_ip^*_{il} 
f^n_{li}\cr
&=\pi_i\sum_{n=1}^{\infty}n^{\gamma}\sum_{m=1}^{\infty}
\sum_{l\neq i} { }_ip^m_{il} { }_ip^n_{li}
=\pi_i\sum_{n=1}^{\infty}n^{\gamma}\sum_{m=1}^{\infty}
{ }_ip^{n+m}_{ii}\cr
&=\pi_i\sum_{n=1}^{\infty}n^{\gamma}\sum_{m>n}
f^{m}_{ii}
=\pi_i\sum_{n=1}^{\infty}\left(\sum_{k=1}^{n-1}k^{\gamma}\right)f^{n}_{ii}\cr
}
$$
and we finish the proof by noting that 
$\sum_{k=1}^{n-1}k^{\gamma} \sim n^{\gamma+1}/ (\gamma +1)$ as
$n\to \infty$.
\hfill$\diamondsuit$
\vsni
\noindent
{\bf Remark.}
It is well known that, for a recurrent chain,
if $M_{ii}^{(\gamma+1)}<\infty$ for some state
$i$ then $M_{ij}^{(\gamma+1)}<\infty$,
for all pairs (distinct or not) $i,j\in S$ (see, e.g., [Chu], Chap. I.11, Cor. 1).
Notice however that even though ${ M}_{\pi i}^{(\gamma)}<\infty$ for all $i\in S$,
the series $\sum_{i\in S}\pi_i { M}_{\pi i}^{(\gamma)}$ is
divergent. To see this, consider for example $\gamma =1$. Assuming 
$M_{ii}^{(2)}<\infty$ let us suppose that 
$$
\sum_{i\in S}\pi_i { M}_{\pi i}^{(1)}=
\sum_{i}\pi_i \sum_{l\not= i}\pi_l M_{li}^{(1)}<\infty
$$
Then, since the double series has positive terms we would have
$$
\sum_{l}\pi_l \sum_{i\not= l}\pi_i M_{li}^{(1)}<\infty,
$$
as well. But this is impossible because 
$M_{il}^{(1)}+M_{li}^{(1)}=(1+{ }_lp^{*}_{ii})\pi_i^{-1}$ ([Chu], p.65) and
$\lim_{i\to \infty}{(M_{il}^{(1)}/ M_{li}^{(1)})}=0$ 
for all $l\in S$ ([Chu], Chap. I.11, Thm. 6; see also [H1]).

We now state the following definition.
\vsni
\noindent
{\bf Definition 1.} {\it Given a recurrent Markov 
chain $P$ with state space $S$, the {\rm ergodic degree} of $P$ is the number}
$$
d=\inf\{\gamma: M_{ii}^{(\gamma+1)}=\infty\;\hbox{\rm for}\;\hbox{\rm some}\; ( \hbox{\rm and}
\;\hbox{\rm then}\;
\hbox{\rm for}\;\hbox{\rm all} )\; i\in S\; \}$$

\noindent
Notice that $M_{ii}^{(0)}=1$ so that the degree satisfies $d> -1$.
In the following we shall 
refer to an ergodic chain as a chain for which $d$
is strictly positive.
If $M_{ii}^{(\gamma)}<\infty$
for every $\gamma$, one says that $P$ has infinite ergodic degree.
This happens for instance if the coefficients $f_{ii}^n$ decay
geometrically with $n$. 
In this case the corresponding 
chain is accordingly called geometrically ergodic. 
We refer to [FMM] for related convergence results 
in the geometrically ergodic case.

\noindent
The preceeding 
observations and Lemma 1 motivate the next definition.
\vsni
\noindent
{\bf Definition 2.} {\it
Given an ergodic chain $P$ with state space $S$ 
and a signed distribution
$\nu$ on $S$, the {\rm $P$-order} of $\nu$ is the number}
$$\sup \{\gamma > 0\, : \, { M}_{\nu i}^{(\gamma)}<\infty
\;\hbox{\rm for}\;\hbox{\rm some}\; ( \hbox{\rm and}
\;\hbox{\rm then}\;
\hbox{\rm for}\;\hbox{\rm all} )\; i\in S\; \}
$$
\vsni
\noindent
{\bf Remark.}
Lemma 1 implies that the ergodic degree of an ergodic chain $P$ 
coincides with the $P$-order of its stationary distribution $\pi$. 
\vsni
\noindent
{\sl Notations:}  Here and in the sequel, 
for two sequences $a_n$ and $b_n$ we shall write
$a_n \sim b_n$
if the quotient $a_n/b_n$ tends to unity as $n\to \infty$. 
Moreover, the notation $a_n={O}_{\epsilon}(n^{-d})$
means that $a_n=o(n^{-(d-\epsilon)})$, $\forall \epsilon >0$,
or, which is the same, that $a_n\cdot n^d$ grows slower than
any power of $n$ as $n\to \infty$. This condition is satisfied 
if, for example, $a_n$ decays as $C\,n^{-d}\,L(n)$
where $L(n)$ is some function slowly varying at infinity, 
i.e. $L(cn) \sim L(n)$ for every positive $c$.
\vsni
\noindent
We now state the main
result of this Section.
\vsni
\noindent
{\bf Theorem 1.} {\sl Suppose $P$ has ergodic degree 
$d > 0$. 
Then, for any initial signed distribution
${\bf \nu}$ of $P$-order at least $d$, we have 
$$
||{\bf \nu} P^n-{\bf \pi}|| ={O}_{\epsilon}(n^{-d}).
$$
In addition, if $M_{ii}^{(d+1)}=\infty$ for some (and then for all) $i$
and the $P$-order of $\nu$ is strictly
larger than $d$,
then the above bound is sharp, i.e.
$n^d\cdot ||{\bf \nu} P^n-{\bf \pi}||$
varies slower than any power of $n$.}
\vsni
\noindent
We let $\tau$ be the shift transformation on $\O$, that is 
$x_k\circ \tau (\o)=\o_{k+1}$. 
With $P$ and ${\bf \pi}$ one can
define a $\tau$-invariant Markov random field $\mu = \mu(P,{\bf \pi})$ 
supported by $\O$ as follows:
$$
\mu(\{ x_k(\o)=\xi_0,\dots ,x_{k+n}(\o)=\xi_{n}\})=
\pi_{\xi_0}\prod_{j=1}^np_{\xi_{j-1}\xi_j}\eqno(1.8)
$$
We shall say that $\mu$ has ergodic degree $d$ whenever 
$P$ (and $\pi$) has the same property.
We then have the following,
\vsni
\noindent
{\bf Corollary 1.} {\sl  Suppose $\mu$ has ergodic degree 
$d > 0$. Then, for any pair of
bounded vectors ${\bf u}, {\bf v} : S \to \R$,}
$$
|\, \mu({\bf u}(x_n){\bf v}(x_0))-\mu({\bf u}(x_0))\,\mu({\bf v}(x_0))\,|
= \, {O}_{\epsilon}(n^{-d})
$$
\vsni
\noindent
{\bf 2. PROOFS}
\vsni
\noindent
We shall prove 
Theorem 1 and its Corollary through several Lemmas.
We start with few technical results which will be used
several times in the sequel.
\vsni
\noindent
{\bf Lemma A.} (see, e.g., [Chu], Chap. I.5)  {\sl Let $\{a_n\}_{n\geq 0}$ be a 
sequence of nonnegative numbers not all vanishing and such that 
$a_n/\left(\sum_{m=0}^na_m\right)\to 0$, $n\to \infty$.
Then, whenever the sequence $\{b_n\}_{n\geq 0}$ of real
numbers has a limit, we have}
$$
\lim_{n\to \infty}{\sum_{m=0}^na_mb_{n-m}\over \sum_{m=0}^na_m}
=\lim_{n\to \infty} b_n.
$$
{\bf Lemma B.}  {\sl Let $D(z)=\sum_{n=0}^\infty d_n z^n$ be
absolutely convergent and $D(z)\neq 0$ for  
$|z|\leq 1$. Let moreover $d_n = {O}_{\epsilon}(n^{-\gamma})$
for some $\gamma \geq 1$. Then 
$$
C (z)={1\over D(z)} = \sum_{n=0}^\infty c_n z^n
$$
is also absolutely convergent for 
$|z|\leq 1$ and $c_n = {O}_{\epsilon}(n^{-\gamma})$.
The assertion remains valid if ${O}_{\epsilon}(n^{-\gamma})$ is replaced 
by $o(n^{-\gamma})$. 

\noindent
If, in addition, $d_0=1$, $d_n> 0$ and $d_n/d_{n-1}$ is increasing,
then $c_0=1$ and $\sum_{k=0}^nc_k>0$ decreases monotonically
to $D(1)^{-1}<1$. }
\vsni
\noindent
{\sl Proof.} 
The first statement is a consequence of a theorem of Wiener and
its proof can be found in [Ro], Lemma 3.II.  For the last statement
see, e.g., [H2], Thm. 22. $\diamondsuit$
\vsni
\noindent
{\bf Lemma C.}  {\sl Let $D(z)$ and $C (z)$ be as 
in the first part of Lemma B with $d_n = {O}_{\epsilon}(n^{-\gamma})$
for some $\gamma > 1$. Assume furthermore that $d_n \geq 0$ and 
$\sum n^{\gamma-1}\, d_n =\infty$.
Given a sequence $e_n$, $n\geq 0$, let  
$h_n=\sum_{k=0}^nc_k\, e_{n-k}$, or else
$$
\sum_{n=0}^\infty h_n z^n ={\sum_{n=0}^\infty e_n z^n\over \sum_{n=0}^\infty d_n z^n}\cdot
$$
\item{(a)} If $e_n = {O}_{\epsilon}(\sum_{\ell >n}d_\ell)$,
then $h_n ={O}_{\epsilon}(\sum_{\ell >n}d_\ell)$, and
the assertion remains valid if ${O}_{\epsilon}(\sum_{\ell >n}d_\ell)$ is replaced 
by $o(\sum_{\ell >n}d_\ell)$.
\item{(b)} If, in addition, $e_n = \sum_{\ell >n}d_\ell$ then $h_n 
\sim D(1)^{-1}\, e_n.$

}
\vsni
\noindent
{\sl Proof.} First,
since $d_n = {O}_{\epsilon}(n^{-\gamma})$ we have from Lemma B that
$c_n ={O}_{\epsilon}(n^{-\gamma})$ as well. 
To show (a) we then notice that
$$
|h_n| \leq \max_{0\leq k\leq n/2}|e_{n-k}|\sum_{k=0}^{n/2}|c_k|+
\max_{n/2\leq k\leq n}|c_k|\sum_{k=0}^{n/2}|e_k|\, .
$$
Therefore, if $\sum |e_k|<\infty$ then
$h_n= {O}(\max\{|c_n|,|e_n|\})$, otherwise $h_n = {O}(e_n)$. 
Indeed, the condition $\sum |e_k|=\infty$ entails
$1<\gamma \leq 2$. Since
$d_n = {O}_{\epsilon}(n^{-\gamma})$ and
$e_n ={O}_{\epsilon}(\sum_{\ell >n}d_\ell)={O}_{\epsilon}(n^{-\gamma+1})$, Lemma B implies that the last term 
in the r.h.s. of the above expression is 
${O}_{\epsilon}(n^{-2(\gamma -1)}) = o(\sum_{\ell >n}d_\ell)$.

\noindent
Let us now prove assertion (b). 
Under the assumption stated there,
Lemma A yields $h_n - e_n\, \sum_{k=0}^nc_k = o(\sum_{k=0}^nc_k)=o(1)$.
But we can say more. The conditions $d_n \geq 0$ and 
$\sum n^{\gamma-1}\, d_n =\infty$ imply that $n^{\gamma-1} \cdot e_n$ decays slower than any
inverse power of $n$. Moreover, let us note that since $e_n = \sum_{\ell >n}d_\ell$ we have
$D(1)-\sum_{n=0}^\infty d_n z^n=(1-z)\sum_{n=0}^\infty e_n z^n$. We then write
$$
\sum_{n=0}^\infty h_n z^n =D(1)^{-1}\sum_{n=0}^\infty e_n z^n  +
{(1-z)\left(\sum_{n=0}^\infty e_n z^n\right)^2\over D(1)\sum_{n=0}^\infty d_n z^n}.
$$
The proof of (b) then reduces to show that the coefficients of the last power
series are $o(e_n)$. To this end we use the following easily checked fact:
$$
\sum_{k=0}^n k^{1-\gamma}(n-k)^{1-\gamma}=\cases{O(n^{3-2\gamma }), &for $1<\gamma <2$,\cr
                                               O(\log n /n), &for $\gamma =2$,\cr
                                                O(n^{1-\gamma}), &for $\gamma >2$.\cr}
$$
By the above, the coefficients of the power series 
$(1-z)\left(\sum_{n=0}^\infty e_n z^n\right)^2$
are ${O}_{\epsilon}(n^{-2(\gamma -1)})$ 
for $1<\gamma \leq 2$ and
${O}_{\epsilon}(n^{-\gamma})$ for $2<\gamma$, therefore $o(e_n)$ in both cases.
The claim now follows by applying the same reasoning as in the proof of (a) to the coefficients
of $(1-z)\left(\sum_{n=0}^\infty e_n z^n\right)^2\sum_{n=0}^\infty c_n z^n$.
 $\diamondsuit$
\vsni
\noindent 
In the following Lemma we shall establish an asymptotic
equivalence which determines the speed of convergence of
the diagonal transition probabilities $p_{ii}^n$ to the stationary
distribution $\pi_i$ in terms of
the ${\bf P}_i$-distribution of the first return time $r^{(i)}_1$.
This will be prove useful to obtain sharp bounds under
appropriate conditions.
\vsni
\noindent
{\bf Lemma 2.} {\sl For a (finitely) ergodic chain $P$ with state space $S$
and
stationary distribution
$\pi$, we have, for any $i\in S$,}
$$
p_{ii}^n-\pi_i\, \sim \,{1\over m_i^2}\,
\sum_{\ell >n}{\bf P}_i\{r_1^{(i)}(\o)>\ell\}.
$$ 
{\sl Proof.} 
We introduce the generating functions
$$
P_{ij}(z)=\sum_{n=0}^{\infty}p^n_{ij}\, z^n,\qquad
F_{ij}(z)=\sum_{n=0}^{\infty}f^n_{ij}\, z^n\eqno(2.1)
$$
and from (1.2) we get the relations (we set $f^n_{ij}=0$ for $n=0$)
$$
P_{ii}(z)={1\over 1-F_{ii}(z)}, \quad 
P_{ij}(z)=F_{ij}(z)P_{jj}(z),\quad i\neq j.\eqno(2.2)
$$
We first show that the function $P_{ii}(z)$ is
analytic in $|z|<1$ and converges at every point of the
unit circle besides $z=1$.
Indeed, recurrence of the state $i$ implies $F_{ii}(1)=1$,
so that $|F_{ii}(z)|<1$ for $|z|<1$ because $f^n_{ii}\geq 0$.
Moreover, $|F_{ii}(z)|<1$ also for $|z|=1$, $z\neq 1$.
This follows from the fact that, since the chain is
aperiodic, ${\rm g.c.d.}\{n, f^n_{ii}\neq 0\}=1$.
Now set
$$
D_{ii}(z)=\sum_{n=0}^{\infty}d^{(n)}_{ii} z^n,  \;\;\;
d^{(n)}_{ii} :=\sum_{k>n}f^k_{ii}={\bf P}_i\{r_1^{(i)}(\o)>n\}\eqno(2.3)
$$
and notice that $D_{ii}(z)$ converges absolutely in $|z|\leq 1$ and
has no zeros on $|z|=1$. In addition
$\sum_{n=0}^{\infty}d^{(n)}_{ii}=m_i$. It then follows from Lemma B
that the function
$$
{1\over D_{ii}(z)}=(1-z)P_{ii}(z)\eqno(2.4)
$$
has a power series expansion which converges absolutely in the closed
unit disk and, moreover, its value at $z=1$ is
$m_i^{-1}=\pi_i$. Set 
$$
{1\over D_{ii}(z)}=:\sum_{n=0}^\infty c^{(n)}_{ii}\,z^n, \qquad
\sum_{n=0}^\infty |c^{(n)}_{ii}|<\infty.\eqno(2.5)
$$
We now observe that the ergodicity assumption
implies that $d^{(n)}_{ii}={\bf P}_i\{r_1^{(i)}(\o)>n\}=o(n^{-1})$. 
We may then
use again Lemma B to obtain 
$c^{(n)}_{ii} =o(n^{-1})$ as well.
By an Abelian theorem (see, e.g., [Chu], p.55 ) we then have
$$
p^n_{ii} =\sum_{k=0}^nc^{(k)}_{ii}\to  \pi_i,\qquad n\to \infty.\eqno(2.6)
$$ 
To obtain more information, we first observe that $m_i\, p^n_{ii}-1$ 
is the coefficient
of $z^n$ in 
$$
H_{ii}(z):=m_i\, P_{ii}(z)-{1\over 1-z} = {E_{ii}(z)\over D_{ii}(z)}
\eqno(2.7)
$$
where 
$$
E_{ii}(z)=\sum_{n=0}^{\infty}e^{(n)}_{ii} z^n, \;\;\; e^{(n)}_{ii} 
:=\sum_{\ell>n}d_{ii}^{(\ell)}=\sum_{\ell>n}
{\bf P}_i\{r_1^{(i)}(\o)>\ell\}.\eqno(2.8)
$$
Now, if the ergodic degree $d$ is finite the conditions of Lemma C-(b) are satisfied for the
sequences $d^{(n)}_{ii}$, $c^{(n)}_{ii}$ and $e^{(n)}_{ii}$. Whence we conclude that
$$
m_i^2\,(p_{ii}^n-\pi_i)\, \sim \,
\sum_{\ell >n}{\bf P}_i\{r_1^{(i)}(\o)>\ell\}.\eqno(2.9)
$$
This finishes the proof. $\diamondsuit$

\noindent
\vsni
\noindent
{\bf Lemma 3.} {\sl Suppose $M_{ii}^{(\gamma)}<\infty$ for some
(and hence for all) $i\in S$ and
for some $\gamma \geq 1$. Then, }
$$
||{\bf \delta}_iP^n-{\bf  \pi}||=o(n^{-(\gamma-1)}).
$$
\noindent
{\sl Proof.} We start noticing that the assumption 
$M_{ii}^{(\gamma)}<\infty$ implies 
${\bf P}_i\{r_1^{(i)}(\o)>n\} = o(n^{-\gamma})$ and
therefore, by Lemma 2, we have 
$$|p_{ii}^n-\pi_i|=o(n^{-(\gamma-1)}).\eqno(2.10)
$$
More generally, it follows from (1.2),
$\sum_{n}f^n_{ij}=1$ and and Lemma A that
$p_{ij}^n\to\pi_j$ as $n\to \infty$.  
Furthermore, as already remarked, 
the condition $M_{ii}^{(\gamma)}<\infty$
implies that $\sum_{n=1}^{\infty} n^{\gamma}f^n_{ij}<\infty$,
for all pairs (distinct or not) $i,j\in S$. 
This and Lemma A, along with the inequality
$$
|p_{ij}^n-\pi_j|\leq \sum_{k=1}^nf_{ij}^k\,|p_{jj}^{n-k}-\pi_j|
+\pi_j\, \sum_{k>n}f_{ij}^k,\eqno(2.11)
$$
imply that the rate of convergence to zero of $|p_{ij}^n-\pi_j|$
is the same as in (2.10).
These properties entail that $P^n$ tends to the matrix
whose rows are $(\pi_1,\pi_2,\dots)$.
To finish the proof we proceed as follows. Having fixed a state $k\in S$
we use (1.1) along with standard decomposition formulae (see [Chu], Chap. I.9) 
to write $p_{ij}^n-\pi_j$ as 
$$
p_{ij}^n-\pi_j=\sum_{m=1}^{n-1} { }_kp_{kj}^{n-m}\, 
(p_{ik}^m-\pi_k) + { }_kp^n_{ij}- \pi_k \, \sum_{m=n}^\infty { }_kp^m_{kj}=I+II+III.
$$
Recalling that $\sum_j{}_kp_{lj}^n=\sum_{m\geq n}f_{lk }^m$ and summing over
$j\in S$ we immediately obtain $\sum_{j\in S}|II|=o(n^{-\gamma})$ and 
$\sum_{j\in S}|III|=o(n^{-(\gamma-1)})$. For the first term we have
$$
\sum_{j\in S}|I|\leq \sum_{m=1}^{n-1}|p_{ik}^m-\pi_k|\sum_{r\geq n-m}f_{kk }^r,
$$
Let us multiply both sides of the above inequality by
$n^{(\gamma-1)}$. Using the fact that
$n\le m(n+1-m)$ if $1\le m\le n$ we get
$$n^{\gamma-1}\sum_{m=1}^{n-1}|p_{ik}^{m}-\pi_{k}|
\sum_{r\ge n-m}f^{r}_{kk}
\le \sum_{m=1}^{n-1}|p_{ik}^{m}-\pi_{k}|m^{\gamma-1}
\sum_{r\ge n-m}f^{r}_{kk}(n+1-m)^{\gamma-1}.$$
Since
$\lim_{p\to\infty}(p+1)^{\gamma-1}\sum_{r\ge p}f^{r}_{kk}=0$
and $\lim_{m\to \infty}|p_{ik}^{m}-\pi_{k}|m^{\gamma-1}=0,$
from Lemma $A$ it follows that
the r.h.s. tends to zero as $n\to \infty$ and therefore
$\sum\limits_{j\in S}|I|=o(n^{-(\gamma-1)}).$
We have thus found that
$$
\sum_{j\in S}|p_{ij}^n-\pi_j| = o(n^{-(\gamma-1)})
$$
and the proof of Lemma 3 is complete.
$\diamondsuit$
\vsni
\noindent
{\bf Lemma 4.} {\sl For any initial signed distribution
${\bf \nu}$ such that ${ M}_{\nu i}^{(\gamma-1)}<\infty$ 
for some
(and hence for all) $i\in S$ (and $\gamma \geq 1$) and
under the hypotheses of Lemma 3, we have}
$$
||\nu P^n - \pi || = o(n^{-(\gamma-1)}).
$$
{\sl Proof.} 
Putting ${\bf \nu} = \sum \nu_l \, \delta_l$ and using the fact
that ${\bf \nu}$ is normalized, i.e. $\sum \nu_l=1$, we write
$$\eqalign{
{\nu} P^n-{\pi}=\sum_l &\, \nu_l \,\, ({\delta}_iP^n-{\pi})+
\sum_{l \neq i}\nu_l \, ({\delta}_lP^n-{\delta}_iP^n)\cr
&=({\delta}_iP^n-{\pi})+
\sum_{l \neq i}\nu_l \, ({\delta}_lP^n-{\delta}_iP^n).\cr }\eqno(2.15)
$$
The $\ell_1$-norm of the first term in the r.h.s. 
is then estimated by
Lemma 3. For the second term we have
$\Vert {\delta}_lP^n-{\delta}_iP^n \Vert = 
\sum_j | p_{lj}^n-p_{ij}^n|$. Using the decompositions
$p_{lj}^n={}_ip_{lj}^n +\sum_{k=1}^{n-1}f_{li}^kp_{ij}^{n-k}$
([Chu], Chap. I.9, Thm. 1) and $p_{ij}^n=\sum_{k=1}^{n-1}f_{li}^k\, p_{ij}^{n}+
\sum_{k=n}^{\infty}f_{li}^k\, p_{ij}^{n}$,
and noting that $\sum_j{}_ip_{lj}^n=\sum_{k\geq n}f_{li}^k$
and $\sum_jp_{ij}^{n}=1$, we obtain
$$
\Vert {\delta}_lP^n-{\delta}_iP^n \Vert
\leq 
2\sum_{k= n}^{\infty}f_{li}^k+
\sum_{k=1}^{n-1}f_{li}^k\sum_j|p_{ij}^{n-k}-p_{ij}^n|.
$$ 
Thus, by (1.6), the norm of the last term in the r.h.s. 
of (2.15) is bounded by
$$
2\sum_{l \neq i}\nu_l \,\sum_{k\geq n}f_{li}^k+
\sum_{k=1}^{n-1}{\bf P}_{\nu}\{r_0^{(i)}=k\}\sum_j
|p_{ij}^{n-k}-p_{ij}^n|.
$$
The assumption that 
${ M}_{\nu i}^{(\gamma-1)}<\infty$ immediately
implies that the first term in the above expression is $o(n^{-(\gamma-1)})$.
As far as the second term is concerned, we may use the inequality
$$
\sum_{j}|p_{ij}^{n-k}-p_{ij}^{n}|\le
\sum_{j}|p_{ij}^{n-k}-\pi_{j}|+
\sum_{j}|p_{ij}^{n}-\pi_{j}|,
$$
and it will suffice to estimate the expression 
$$
\sum_{k=1}^{n-1}{\bf P}_{\nu}\{r_0^{(i)}=k\}\sum_j
|p_{ij}^{n-k}-\pi_{j}|.
$$
Now, the assumption 
${ M}_{\nu i}^{(\gamma-1)}<\infty$ implies that 
$\lim_{k\to\infty}k^{\gamma-1}{\bf P}_{\nu}\{r_0^{(i)}=k\}=0$ and,
 under the assumptions of Lemma 3,
$\lim_{m\to\infty}m^{\gamma-1} \sum_j
|p_{ij}^{m}-\pi_{j}|=0$.
We may then repeat the argument given at the end of the proof of Lemma 3 to see
that the above expression is $o(n^{-(\gamma-1)})$.
\hfill $\diamondsuit$
\vsni
\noindent
{\sl Proof of Theorem 1.}
The conditions on the ergodic degree of $P$ and on the
$P$-order of $\nu$ imply that the assumptions of Lemmas 3 and 4
are satisfied for $\gamma = d+1-\epsilon$, $\forall \epsilon >0$.
This gives a rate of convergence $o(n^{-(d-\epsilon)})$,
$\forall \epsilon >0$, that is ${O}_{\epsilon}(n^{-d})$.
But we can say more. Indeed, the condition 
$\sum n^{d+1}f_{ii}^n=\infty$ and Lemma 2 entail that
$\vert p_{ii}^n-\pi_i\vert \cdot n^d$, and thus 
$||\delta_iP^n-\pi||\cdot n^d$, 
decays
slower than any inverse power of $n$.
On the other hand, from the proof of Lemma 4 we see that the condition
that $\nu$ has $P$-order strictly 
larger than $d$ implies that the norm of
$\sum_{l \neq i}\nu_l \, ( {\delta}_lP^n-{\delta}_iP^n)$
is ${O}_{\epsilon}(n^{-d'})$, for some $d'>d$. This
prevents from possible cancellations
among the two terms in the r.h.s. of (2.15). 
\hfill$\diamondsuit$
\vsni
\noindent
{\bf Remark.} The proof given 
above brings out the meaning
of the condition on the $P$-order of the initial distribution
$\nu$. This is related 
to the fact that the behaviour of
$|p_{ij}^n-\pi_j|$ and hence of
$||{\bf \delta}_iP^n-{\bf  \pi}||$ is 
necessarily not uniform in the
departing state index $i$. Indeed, according to the
above discussion, such uniformity would
imply the existence of two positive constant $C_1, C_2$
and an integer $n_0$,
which {\sl do not depend} on $i$ and $l$, 
such that, for all $n\geq n_0$ 
$$
C_1 \leq {\sum_{k\geq n}f_{li}^k \over 
\sum_{k\geq n}f_{il}^k} \leq C_2.
$$
This, in turn, would imply that the ratio 
$M_{li}^{(1)}/M_{il}^{(1)}$ satisfies a similar bound. 
On the other hand, as already observed,
$\lim_{l\to \infty}{(M_{li}^{(1)}/ M_{il}^{(1)})}=0$, 
for all $i\in S$.
\vsni
\noindent
{\sl Proof of Corollary 1.}
For any pair ${\bf u}\in \ell_{\infty}(S)$, ${\bf \rho} \in \ell_1(S)$ we define 
${\overline {\rho {\bf u }}}=
(\rho(1)u(1),\rho(2)u(2),\dots)$ and
$ \rho\cdot{\bf  u} = \sum_{i\in S}
\rho(i)u(i)$. Thus ${\overline {\rho{\bf  u}}}\cdot {\bf 1} =
{\bf \rho} 
\cdot {\bf u}$,
and the unit column vector ${\bf 1}=(1,1,\dots)^t$ satisfies 
$P{\bf 1}={\bf 1}$. For definiteness and without loss, 
suppose that $\mu({\bf u})\, \mu({\bf v})\neq 0$.
Then we have
$$
\eqalign{
|\, \mu({\bf u}(x_n){\bf v}(x_0))-
\mu({\bf u}(x_0))\,\mu({\bf v}(x_0))\,| &=
|\,{\overline {{\bf \pi} {\bf v}} } \, P^n\cdot {\bf u} - 
({\pi} \cdot {\bf v})({\pi} \cdot {\bf u}) \, | \cr
&=|\,({\overline {{\bf \pi} {\bf v}}}\, P^n-
{\bf \pi} \, ({\overline {{\bf \pi} {\bf v}}}\cdot {\bf 1})\,)
\cdot {\bf u}\, |\cr  
&\leq \Vert {\bf u}\Vert_{\infty}\, \Vert {\bf v}\Vert_{\infty}\,
\Vert \,\nu\,P^n-{\bf \pi}\, \Vert\cr 
}
$$
where $\nu$ denotes the normalized $\ell_1$ row vector 
${\overline {{\pi} {\bf v}}}/ ({\pi} \cdot {\bf v})$.
The result now follows putting together Lemma 1 and Theorem 1.
\hfill$\diamondsuit$
\vskip 0.5cm
\noindent
{\bf 3. CONVERGENCE VS ANALYTIC AND SPECTRAL PROPERTIES. 
AN EXAMPLE.}
\vsni
\noindent
As we have seen, the dependence 
on the departing state $i$ of the behaviour of
$\Vert \delta_i P^n -\pi \Vert$, 
although not explicitly indicated in Lemma 3,
is what makes
our assumptions on the $P$-order of the initial 
distribution $\nu$ necessary.

\noindent 
Moreover, from our discussion it follows that 
the rate of convergence
to zero of $\Vert \delta_i P^n -\pi \Vert$ is connected
with the analytic properties of the generating functions 
$P_{ij}(z)$ in the vicinity of the singular point $z=1$. 

\noindent
If we now consider $P$ as a bounded linear 
Markov operator acting on the Banach space $\ell_1(S)$, 
its adjoint $P^*$ is
represented by the transposed matrix acting on the dual space 
$\ell_1^*=\ell_{\infty}$. The resolvent
$R_{\lambda}(P):=(\lambda I - P)^{-1}$ admits, 
for $|\lambda|> \Vert P\Vert$, the expansion
$$
\lambda \, R_{\lambda}(P) = I + 
\sum_{n=1}^{\infty}\left({P\over \lambda}\right)^n
$$
which shows that $1-\delta_{ij}+P_{ij}(z)$ is the 
$(i,j)$-element of $\lambda \, R_{\lambda}(P)$, 
with the identification $z=1/\lambda$. 
This, in turn, indicates that the convergence properties
of $\Vert \delta_i P^n -\pi \Vert$, the analytic properties
of the functions $P_{ij}(z)$, and the spectral properties
of $P$ in $\ell_1(S)$ are intimately connected items.
In particular, the dependence of the first two from the 
state index $i$ plays an important role in 
determining nature of the latter,
as we shall see in the following example\footnote{$^{1}$}{
We shall adopt the convention that a matrix
$(t_{ij})$ representing an operator $T$ acts from the right, that is
through the equations $(Tx)_j=\sum_{i\in S}x_i\,t_{ij}$.}.

\vsni
\noindent
{\sl Example.} 
Suppose that $S=\N$ and the transition matrix is
$$
P=\pmatrix{p_1&p_2&p_3&\ldots \cr
           1  &0  &0  &\ldots \cr
           0  &1  &0  &\ldots \cr
           0  &0  &1  &\ldots \cr
           \vdots &\vdots &\vdots&\ddots \cr}
$$
The space $\O$ is then given by all sequences $\o$ satisfying
the following condition: given $\o_i$ then either $\o_{i-1}=\o_i+1$
or $\o_{i-1}=1$.
We shall assume that the probability vector $p=(p_1,p_2,\dots)$
has the property
${\rm g.c.d.}\{n: p_n>0\}=1$.
It then follows that the corresponding 
chain is aperiodic and recurrent.
Let the coefficients $d_n$ be defined by
$d_n:=\sum_{i>n}p_i$, ($n\geq 0$). 
The steady-state equation is $\pi_n=\sum_{i\in S}\pi_i\, p_{in}$
and is formally solved by $\pi_n=\pi_1\, d_{n-1}$, ($n\geq 1$).
We also have $f^n_{11}=p_n$.
Consequently, the chain is positive-recurrent
if and only if $\sum d_n < \infty$, 
null-recurrent in the opposite case. 
In the former case, we have
$\pi_1=(\sum_{n=1}^{\infty}np_n)^{-1}=
(\sum_{n=0}^{\infty}d_n)^{-1}$.
Notice that the two probability vectors $\pi$ and $p$ coincide
if and only if $p_n=2^{-n}$. On the other hand, if 
$p_n\sim n^{-(d +2)}\, L(n)$ with $L(n)$ a suitable function 
slowly varying at infinity 
then the chain
has ergodic degree $d$.

\vsni
\noindent
{\bf Remark 1.} It is not difficult to realize that
the $\tau$-invariant Markov random field $\mu = \mu(P,{\bf \pi})$ 
defined in (1.8), with $P$ and ${\bf \pi}$ as above, can
be viewed as an equilibrium state [Ru] for the 
continuous potential function
$V: \O \to \R$ defined as
$$
V(\o) =\log p_{\o_0} - \log p_{\o_1} + \log P(\o_0,\o_1).
$$
\vsni
\noindent
{\bf Remark 2.} The Markov chain $P$ is a 
reference model in renewal theory (see [Se]). In particular,
the validity of the renewal limit theorem corresponds to
the fact that the chain is ergodic. Several estimates
on the remainder term in this limit theorem (which corresponds
to the speed of convergence to equilibrium) have been
obtained. See [Ro] for very accurate results 
and also [Se], Chap. 24, for a review.
These results can be viewed as particular cases
(corresponding to $\nu=\delta_i$ and $u_k=\delta_i^k$,
for some $i\in \N$)
of Theorem 2.III stated below.
Moreover, this example has interesting applications 
in modelling renewal processes arising in dynamical system
theory; a situation which has recently become a standard example
being that of Markov
interval maps modelling
temporal intermittency (see, e.g., [Wa]). A brief discussion
on the consequences of the results stated below 
in the context of dynamical systems theory
is given in the Appendix at the end of the paper.

\vsni
\noindent
{\bf Theorem 2.} {\sl Suppose that the chain $P$ defined above
has finite ergodic degree $d>0$. Then,
\item{I.} The generating functions $P_{ij}(z)$ 
defined in (2.1) are analytic in
the open unit disk. For $|z|\leq 1$ the functions $1/P_{ij}(z)$ have only one zero at $z=1$
which is a non-polar singular point for $P_{ij}(z)$.
\item{II.} The spectrum $\sigma (P)$ of 
the Markov operator $P$ acting on 
$\ell_1(\N)$ coincides with the closed unit disk and 
decomposes as follows:
$\sigma_p (P)=\{\lambda : |\lambda|<1\}\cup \{1\}$ and
$\sigma_c (P)=\{\lambda : |\lambda|=1, \lambda \neq 1\}$.
\item{III.} For any 
bounded vector ${\bf u}$ and any initial
distribution $\nu=(\nu_i)_1^{\infty} \in \ell_1(S)$ s.t. 
$\nu_i = { O}(\pi_i)$, 
the quantity 
$({\bf \nu} P^n-{\bf \pi})\cdot {\bf u}$
decays as ${O}_{\epsilon}(n^{-d})$. 

\noindent
Assume furthermore that 
$p_n\sim n^{-(d +2)}\, L(n)$ with $L(n)$ 
slowly varying at infinity and $u_i = o (1)$, $\nu_i=o (\pi_i)$. 
Then we have
$$
({\bf \nu} P^n-{\bf \pi})\cdot {\bf u}\sim 
C \, n^{-d}\, L(n),
$$
with $C=(\pi \cdot {\bf u})(\nu \cdot {\bf 1})/(d(d+1)m_1)$.}

\vsni
\noindent
{\bf Remark 1.} 
Statement $I$ above holds for any aperiodic Markov chain with finite
ergodic degree and is well known. On the other hand, it
can be considerably improved
by specifying further properties of the probability vector $p$.
For instance, if the $p_n$ form a monotonically decreasing sequence
 $p_1\geq p_2\geq \cdots$ satisfying the Kaluza property:
$p_n^2>p_{n+1}\, p_{n-1}$ (with $p_0=1$) then using the last part of Lemma B one can show that the generating functions
$P_{ij}(z)$ can be continued meromorphically to the entire
$z$-plane with a branch cut along the ray
$(1,+\infty )$ (see [Is2]).
\vsni
\noindent
{\bf Remark 2.}
In the null-recurrent case ($d\leq 0$)
the statements corresponding to II and III above are modified as follows (see [A]):
{\it \item{II'.} The spectrum $\sigma (P)$ of 
the Markov operator $P$ acting on 
$\ell_1(\N)$ coincides with the closed unit disk and 
decomposes as:
$\sigma_p (P)=\{\lambda : |\lambda|<1\}$,
$\sigma_c (P)=\{\lambda : |\lambda|=1, \lambda \neq 1\}$ and
$\sigma_r (P) =\{1\}$.
\item{III'.} Let ${\bf v}=(v_i)_1^{\infty} \in \ell_{\infty}(S)$ be the unique (non-normalized) positive 
invariant vector for $P$ with $v_1=1$ (see [De], Thm 1). Here
$v_n=d_{n-1}$. For any vector ${\bf u}\in \ell_\infty (S)$ such that
${\bf u}\cdot {\bf v}<\infty$ 
 and any initial
distribution $\nu\in \ell_1(S)$ we have
$$
{\bf \nu} P^n\cdot {\bf u}\sim ({\bf \nu}\cdot {\bf 1})({\bf u}\cdot {\bf v}) \, p_{11}^n
$$
and $p_{11}^n\cdot n^{-d}$ varies slower than any power of $n$.
}

\vsni
\noindent
The proof of Theorem 2 will follow from the points 
$I$, $II$ and $III$ discussed hereafter. 

\vsni
\noindent
{\sl I. Generating functions.} 

\noindent
First, it is easy to check that
all entries of the
first $n$ rows of $P^n$ are positive, the $i$-th row of
$P$ being the $(i+n-1)$-th of $P^n$. More specifically, one sees inductively 
that for $n>1$, $i>1$, $j\in \N$, 
$$
P^{n}(i,j)=P^{n-1}(i-1,j).\eqno(3.1)
$$
For the generating functions 
of the $P^{n}(i,j)$'s we then obtain the relations
$$\eqalign{
P_{ij}(z) &= \delta_{ij}+ z^{i-1}P_{1j}(z),\quad j\geq i>1 \cr
P_{i1}(z) &= z^{i-1}P_{11}(z),\quad i\geq 1 \cr
P_{ij}(z) &= z^{i-j} + z^{i-1}P_{1j}(z),\quad i>j>1.\cr }
\eqno(3.2)
$$
It then suffice to study the behaviour of the 
entries of the first row.
They satisfy the recurrence relations
$P^{n}(1,j)=P^{n-1}(1,1)\, P(1,j) + P^{n-1}(1,j+1)$, $j\geq 1$
(recall that 
$P^0(i,j)=\delta_{ij}$).
This yields
$$
P^{n}(1,j)=\sum_{k=1}^{n-1}P(1,k)\ P^{n-k}(1,j) + P(1,j+n-1).
\eqno(3.3)
$$
Putting $j=1$ and recalling that $P(1,k)=p_k=f^k_{11}$
one gets a particular case of equation (1.2). It hence follows that
$$
P_{11}(z) = {1\over 1-\sum_{n=1}^{\infty}p_nz^n}
={1\over (1-z)D(z)}\eqno(3.4)
$$
where $D(z)=\sum_{n=0}^{\infty}d_nz^n$. More generally,
we get for $j>1$
$$
P_{1j}(z) = {z^{1-j}\, P_j(z)\, P_{11}(z)}\eqno(3.5)
$$
where $P_j(z)=\sum_{n=j}^{\infty}p_{n}z^n$.
Finally, using (1.1)-(1.2) along with (3.2), (3.4) and (3.5) we obtain
$$\eqalign{
F_{ij}(z) &= z^{i-j},\qquad\qquad\qquad \quad i> j, \cr
F_{ij}(z) &={z^{i-j}P_j(z)\over  1-\sum_{ 0<n<j}p_nz^n}, \quad j\geq i. \cr}
\eqno(3.6)
$$
\vsni
\noindent
{\bf Remark.} 
As an application of the above formulas one can compute
the moments $M_{ij}^{(\gamma)}$ of $P$.
For instance, if $d>1$, computing the second derivative at $z=1$ of
$F_{ii}(z)$ yields
$$
M_{ii}^{(2)} = {\pi_1\over \pi_i}\left(M_{11}^{(2)}+ 
{2\sum_{n=1}^{i-1}np_n\over \pi_{i}}\right) \sim {2\over \pi_i^2}
$$
where $M_{11}^{(\gamma)}=\sum n^{\gamma}p_n$ and the last asymptotic equivalence
holds for $i\to \infty$.

\vsni
\noindent 
The proof of the analytic properties of the generating functions $P_{ij}(z)$  now follows
a standard path and we therefore omit it.
\vsni
\noindent
{\sl II. Spectral properties of $P:\ell_1(\N)\to \ell_1(\N)$.} 

\noindent
From (3.1)-(3.2) we have that 
the rate of convergence of $P^n(i,j)$ to $\pi_j$ is not uniform
in the departing state $i$ 
(see also the Remark after the proof of Theorem 1). 
We are now going to see how this fact
reflects in the nature of the spectrum of $P$ in $\ell_1$.
In particular, the eigenvalue $1$ is
not isolated, even in the case where
the $p_n$'s are exponentially decreasing.

\noindent
We study the structure of the
spectrum of $P$ using the method of generating functions
(see, e.g., [VJ]). Setting $x=(x_1,x_2,\dots)$ and
$X(w)=\sum_{n=1}^{\infty}x_nw^n$ the formal solutions to the vector equations
$$
(\lambda I - P)x=0 \quad\hbox{and}\quad (\lambda I - P^*){\bf x}=0
$$
can be written as
$$
X(w)={x_1w(1-w)D(w)\over 1-\lambda w}\eqno(3.7)
$$
and 
$$
X(w)=\, {w\over \lambda -w}\; p\cdot {\bf x},\eqno(3.8)
$$
respectively, where $p\cdot {\bf x} = \sum_{n\geq 1}x_np_n$.
The equation $1-\lambda w=0$ (and its reciprocal $\lambda -w=0$)
entails that the boundary of $\sigma(P)$ (and of $\sigma(P^*)$)
is the unit circle. 
\noindent
Let us first consider the point $\lambda =1$. The formal
expressions in (3.7) and (3.8) become 
$$
X(w) = x_1 w\,D(w) \quad\hbox{and}\quad 
X(w)={w\over 1-w}\; p\cdot {\bf x}.\eqno(3.9)
$$
The latter has the solution $X(w)=w/(1-w)$ which is the generating
function of the unit vector in $\ell_{\infty}$. On the other hand,
the former is the generating function of an $\ell_1$-vector
if and only if $D(1)<\infty$. 
Hence, we have that in the positive-recurrent 
case $1$ lies in $\sigma_p(P)$ (for the null-recurrent chain
it lies in $\sigma_r(P)$). 

\noindent
More generally, from (3.7) and (3.8) one sees that
the open unit disc $\{\lambda : |\lambda |<1\}$
is always in the point spectrum. Indeed, 
the function $(1-w)D(w)=1-\sum_{n=1}^\infty p_n w^n$ appearing 
in (3.7) is absolutely convergent for $|w|\leq 1$. 
If $|\lambda |<1$ the same holds true for 
the function $w/(1-\lambda w)=\sum_{n=1}^\infty\lambda^{n-1} w^n$. 
Therefore 
the power series expansion of $X(w)$, being the product of two absolutely convergent
power series, is absolutely convergent at any point of the closed unit disk
$|w|\leq 1$. More precisely, an easy calculation shows that for $n\geq 2$ the coefficient $x_n$ of $w^n$ is bounded above by
$|x_1|\, (|\lambda|^{n-1} + \sum_{k=0}^{n-2}|\lambda|^k p_{n-k-1})$.
This shows that for any $|\lambda |<1$ the function $X(w)$
is the generating function of a vector $x\in \ell_1$.
A similar reasoning shows that for any $|\lambda |<1$ the function
$X(w)$ in 
(3.8) is the generating function of a vector in $\ell_{\infty}$, 
thus proving that $\{\lambda : |\lambda |<1\}\subseteq \s_p(P)$.

\noindent
We conclude by showing that any $\lambda$ s.t.
$|\lambda|=1$, $\lambda \neq 1$ lies in $\sigma_c(P)$.
Indeed, take $\lambda = e^{i\theta}$ with $0 < \theta < 2\pi$ and assume that 
$(\lambda I -P^*)x=0$
for some $x\in \ell_{\infty}$. 
Then the equation in (3.8) gives for the coefficients $x_n$
the relation $x_n = e^{-i(n+1)\theta}\, p\cdot {\bf x}$. So, if $x\ne 0$, then 
$p\cdot {\bf x}\ne 0$.
Multiplying by $p_n$ and summing over $n$ we then get
$1=\sum_n p_n e^{-i(n+1)\theta}$ which is impossible in our case. If the point $\lambda$ 
belongs to the unit circle and is different from $\lambda =1$, 
then the generating function $X(w)$ in (3.7) tends to infinity as $w$ approaches
$\lambda^{-1}$ because $D(w)\ne 0$ for any $|w|=1$. But if the solution
$x$ to the equation $(\lambda I -P)x=0$ belongs to $\ell_1$, then the generating
function $X(w)$ is absolutely convergent at any point of the unit circle and its absolute
value is bounded by $|x|_1$. We then see that the point $\lambda$ does not belong
neither to $\s_p(P^*)$ nor to $\s_p(P)$. This means that $\lambda \in \s_c(P)$.
In particular, we have found that the eigenvalue $1$ is not isolated but is embedded in a continuous spectrum.

\vsni
\noindent
{\sl III. Convergence properties.} 

\noindent
Next, we discuss the convergence properties of this
chain under the hypothesis that it is positive-recurrent.
Note that the
first part of statement III in Theorem 2 is a consequence of Theorem 1, for
$|({\bf \nu} P^n-{\bf \pi})\cdot {\bf u}|\leq \Vert {\bf \nu} P^n-{\bf \pi}\Vert_1 \cdot 
\Vert{\bf u}\Vert_\infty$.
Nevertheless, we shall give an alternative proof which on the one hand 
yields the actual asymptotic behaviour under the hypotheses stated
in the second part of Theorem 2-III and on the other hand
allows us to introduce
a method which appears to be
interesting in its own, for it may be extended to some
more general (i.e. non-markovian) mixing Gibbs 
random fields [Is1].

\noindent
For $z\in \C$, consider the matrix $L_z$ given by
$$
L_z=\pmatrix{p_1z&p_2z&p_3z&\ldots \cr
           p_1z^2  &p_2z^2  &p_3z^2  &\ldots \cr
           p_1z^3  &p_2z^3  &p_3z^3  &\ldots \cr
           \vdots &\vdots &\vdots \cr}
$$
For $z=1$ the matrix $L_z$ can be viewed as the transition
matrix of the process $r_0^{(1)},r_1^{(1)},\dots$
given by the sequence of times between returns to the
state $1$ (see (1.4)).
The vector equation $y = L_z x$, takes the 
generating function form
$Y(w)=p\cdot {\bf 1}_w\, X(z)$ where $Y(w)=\sum_{n=1}^\infty y_n w^n$ and
${\bf 1}_w=(w,w^2,w^3,\dots)^t$. Therefore
the power series of
$L_z$ when acting on $\ell_1(S)$ converges absolutely for any $z$ in the closed unit disk
$|z|\leq 1$.
In addition, there is a simple algebraic relation between the matrices $L_z$ 
and $P$: let $Q$ be the transient chain given by the matrix
$$
Q=\pmatrix{0  &0  &0  &\ldots \cr
           1  &0  &0  &\ldots \cr
           0  &1  &0  &\ldots \cr
           0  &0  &1  &\ldots \cr
           \vdots &\vdots &\vdots&\ddots \cr}
$$
An easy calculation shows that
$$
(I-zQ)(I-L_z)=(I-zP).\eqno(3.10)
$$
This relation entails that if $u$ is an eigenvector of
$P$ with eigenvalue $1/z$, then 
$v=u(I-zQ)$ is an eigenvector
of $L_z$ with eigenvalue $1$. On the other hand we already know that
$P$, when acting on $\ell_1$, has spectral radius equal to $1$ and
no eigenvalues on the unit circle besides eventually $1$.
The choice $z=1$ gives $u=\pi$ and 
$v={\bf \pi}(I-Q)=\pi_1\,p$, as expected.

\noindent
Let now ${\bf u}:S\to R$ be a bounded vector and ${\bf \nu}$
an initial distribution on $S$,
which will be assumed to decay not slower than $\pi$
at infinity. The latter condition is equivalent to the assumption
made in Theorem 1: if the $P$ has ergodic degree
$d>0$ then $\pi$ ($\nu$) has $P$-order (at least) $d$. 

\noindent
Let us consider the following generating function,
$$
S(z)=\sum_{n=0}^{\infty}z^n \, 
({\bf \nu} P^n-{\bf \pi})\cdot {\bf u}.
$$
Using (3.10) we get for $|z|<1$,
$$
\sum_{n=0}^{\infty}z^n \, {\bf \nu} P^n \cdot {\bf u}= 
{\bf \nu} (I-zP)^{-1}\cdot {\bf u}
={\bf \nu} (I-L_z)^{-1}(I-zQ)^{-1}\cdot {\bf u}\, .
$$
Now observe that ${\bf \nu} \, L_z = ({\bf \nu}\cdot {\bf 1}_z) \, p$. 
Iterating $n$ times
we get 
${\bf \nu} \, L^n_z = ({\bf \nu}\cdot {\bf 1}_z) \, \lambda_z^{n-1}\, p$,
with $\lambda_z= p\cdot {\bf 1}_z$, and
the above expression becomes
$$
{({\bf \nu}\cdot {\bf 1}_z) \, p\, 
(I-zQ)^{-1}\cdot {\bf u}\over 1-\lambda_z} 
+{\bf \nu}\,(I-zQ)^{-1}\cdot {\bf u} 
=
{({\bf \nu}\cdot {\bf 1}_z) \, 
( m_1\,{\bf \pi}_z\cdot {\bf u}) \over 1-\lambda_z} 
+{\bf \nu}_z\cdot {\bf u} 
$$
where $m_1=\pi_1^{-1}=D(1)$, ${\bf \nu}_z={\bf \nu}\,(I-zQ)^{-1}$  and 
${\bf \pi}_z=\pi_1 p\,(I-zQ)^{-1}$  
(in particular ${\bf \pi}_z|_{z=1} \equiv \pi$). Therefore
a short manipulation yields the expression
$$
S(z)=({\bf \pi}\cdot {\bf u})\,({\bf \nu}\cdot {\bf 1})\,
H(z) + R(z)
$$
where 
$$
H(z)= {m_1 \over 1-\lambda_z}-{1\over 1-z}= 
{\sum_{n=0}^{\infty}e_nz^n\over \sum_{n=0}^{\infty}d_nz^n},
\quad\hbox{with}\quad e_n = \sum_{k>n}d_k,
$$
and
$$
R(z)={(m_1{\bf \pi}\cdot {\bf u})\,
({\bf \nu}\cdot {\bf 1}_z -{\bf \nu}\cdot {\bf 1})
+({\bf \nu}\cdot {\bf 1}_z) \, (m_1{\bf \pi}_z\cdot {\bf u}-
m_1{\bf \pi}\cdot {\bf u})\over 1-\lambda_z}
+{\bf \nu}_z\cdot {\bf u}.
$$
Using the above and Lemma C one sees that
if $P$ has ergodic degree $d$ then the coefficients of $H(z)$
decay as $\pi_1\, e_n={O}_{\epsilon}(n^{-d})$. It remains to examine the 
behaviour of $R(z)$. 
We have
$$
{{\bf \nu}\cdot {\bf 1} -{\bf \nu}\cdot {\bf 1}_z
\over 1-\lambda_z}={\sum_{n=0}^{\infty}\eta_nz^n\over \sum_{n=0}^{\infty}d_nz^n},
\quad\hbox{with}\quad\eta_n = \sum_{k>n}\nu_k.
$$
Moreover, a straightforward calculation yields
$$
m_1{\bf \pi}_z\cdot {\bf u}=
\sum_{n=0}^{\infty}z^n\left(\sum_{k=1}^{\infty}u_k\,p_{n+k}\right)
$$
and therefore 
$$
{ m_1{\bf \pi}\cdot {\bf u}-
m_1{\bf \pi}_z\cdot {\bf u}\over 1-\lambda_z}=
{\sum_{n=0}^{\infty}\xi_nz^n\over \sum_{n=0}^{\infty}d_nz^n},
\quad\hbox{with}\quad\xi_n = \sum_{k=1}^{\infty}u_k\, d_{k+n}.
$$
In addition,
$$
{\bf \nu}_z\cdot {\bf u}= 
\sum_{n=0}^{\infty}\gamma_n\, z^n\quad\hbox{with}\quad\gamma_n =\sum_{k=1}^{\infty}u_k\,\nu_{n+k}.
$$ 
On the other hand,
$$
|\xi_n| \leq \Vert u\Vert_{\infty}\sum_{k>n}d_k=
\Vert u\Vert_{\infty}\,e_n,\qquad
|\gamma_n| \leq \Vert u\Vert_{\infty}\sum_{k>n}\nu_k=
\Vert u\Vert_{\infty}\,\eta_n.
$$
Reasoning as in the proof of Lemma C we have that if
$\sum |\xi_n|<\infty$ then the coefficient
of $z^n$ of the product $\sum_{n=0}^{\infty}\xi_nz^n \cdot 
\sum_{n=0}^{\infty}\nu_nz^n$ is 
${O}(\max\{|\xi_n|,|\nu_n|\})$ (recall that $\nu_i = { O}(\pi_i)$), otherwise it is ${O}(\xi_n)$. 
Therefore, by the first estimate above, it is ${ O}(e_n)$ in both cases.

\noindent
Comparing all the terms above and using again Lemma C
we have found that under our assumptions
on the distribution ${\bf \nu}$ and the vector ${\bf u}$, 
the quantity $({\bf \nu} P^n-{\bf \pi})\cdot {\bf u}$
decays as ${O}_{\epsilon}(n^{-d})$. 

\noindent
We conclude by deriving the exact asymptotic behaviour of
$({\bf \nu} P^n-{\bf \pi})\cdot {\bf u}$ under the additional
hypotheses imposed in the last part of Theorem 2.III. 
First, if $p_n \sim n^{-(d+2)}L(n)$ then we have
$d_n \sim (d +1)^{-1}n^{-(d +1)}\, L(n)$ and
$e_n \sim d^{-1}(d +1)^{-1}n^{-d}\, L(n)$.
Lemma C then implies that the coefficients of the power series of
$H(z)$ are asymptotically equivalent to $d^{-1}(d +1)^{-1}D(1)^{-1}n^{-d}\, L(n)$.
Moreover, if
$u_i = o (1)$ and $\nu_i=o (\pi_i)$ then
$\xi_n = o (e_n)$ and $\eta_n = o (e_n)$. Again by virtue of Lemma C
this prevents from possible cancellations among
the various coefficients introduced above and yields 
the claim.
\vfill\eject
\noindent
{\bf APPENDIX.}  $\;$ {\sl Renewal chains and
Markov approximations of dynamical systems.} 
Let $(X,\rho)$ be a probability space and $f:X\to X$
be a transformation preserving the
probability measure $\rho$ which we assume to be ergodic. 
Given a measurable subset
$E\subset X$, the quantity 
$$
e_n = {\rho(E\cap f^{-n}E)\over \rho(E)}
\eqno(A.1)
$$
is the probability to observe a return in $E$ after $n$ 
iterations of $f$ (for the first time or not).
The return time function
$$
R_E(x)=\inf\{n>0 : f^n(x)\in E\}\eqno(A.2)
$$ 
is defined (and finite) for a.e. $x\in E$. $E$ itself becomes
a probability space with measure $\rho_E(A)=\rho(A\cap E)/\rho (E)$.
One may then define the induced transformation
$$
f_E(x)=f^{R_E(x)}(x)\eqno(A.3)
$$ 
for a.e. $x\in E$. Both $R_E$ and $f_E$ are measurable and
in fact it is not difficult 
to check that $f_E$ preserves the measure
$\rho_E$ which is of course ergodic. 
We denote by $E_n=\{x\in E : R_E(x)=n\}$ 
the $n$-th levelset of $R_E$. Notice that the above 
construction yields a countable partition
${\cal A}=\{A_n\}$ of $X$
into the sets 
$$
A_n= f^{-(n-1)}(E)\setminus  (\cup_{k=0}^{n-2}f^{-k}(E))
=\cup_{k\geq n} f^{k-n+1}(E_k)\eqno(A.4)
$$
and, $\rho$ being $f$-invariant,
$$
\rho (A_n) =\sum_{k\geq n}\rho (E_k). \eqno(A.5)
$$
Therefore we have
$1=\rho (X) = \sum \rho (A_n) =\sum n\, \rho(E_n)$. 
It hence follows that
$$
\rho_E(R_E)=1/\rho(E),\eqno(A.6)
$$
which is a version of Kac's formula.
Now notice that the number $e_n$ may be rewritten as
$$
e_n = \rho_E(f^{-n}E).\eqno(A.7)
$$
This expression allows us to give another interpretation
of $e_n$. For $x\in E$, let $S_n(x)=\sum_{k=0}^{n-1}R_E(f_E^k(x))$
be the total number of iterates of $f$ needed to observe $n$
returns to $E$ and $N_n(x)=\sum_{k=1}^{n}\chi_E(f^k(x))$
the number of returns up to the $n$-th iterate of $f$. 
A short reflection gives that $\rho_E(S_k\leq n)=\sum_{r=k}^n\rho_E(N_n=r)$. In addition
we have $\rho_E(S_k=n)=\rho_E(S_k\leq n)-\rho_E(S_k\leq n-1)$ for
$k<n$ and $\rho_E(S_n=n)=\rho_E(S_n\leq n)$. A straightforward
computation using these
observations and (A.7) yields (for $n>0$):
$$e_n= \sum_{k=1}^{n}\rho_E(S_k=n)=
\rho_E(N_n)-\rho_E(N_{n-1}),\eqno(A.8)
$$
where $\rho_E(N_n)$ denotes the mean of the random variable $N_n$
(we set $N_0=0$). Thus, $e_n$ may be regarded as the
expected number of returns in $E$ per iteration of $f$ (after
$n$ iterations). It then turns out that the validity of
the renewal theorem for $e_n$, that is [Se]:
$$
e_n \to {1\over \rho_E(R_E)},\qquad n\to \infty \eqno(A.9)
$$
is equivalent to the (self-)mixing property for the set $E$, that is
$e_n \to \rho (E)$. A further remark is the following. 
Let us decompose
$$\eqalign{
e_n&=\sum_{k=1}^n \rho_E(f^l(x)\notin E, 0<l<k, f^k(x)\in E, f^n(x)\in E)\cr
&=
\sum_{k=1}^n \rho_E(E_k)\cdot 
\rho_E(f^n(x)\in E\, | \, R_E(x)=k)\cr}\eqno(A.10)
$$
Now suppose that the process $\{f^n(x)\}$ ``renews" itself 
each time it returns to $E$. In other words, suppose that the 
random variables $R_E, R_E\circ f_E, R_R \circ f_E^2, \dots$
defined on the probability space $(E,\rho_E)$ are mutally
independent. In this case we would have
$$
\rho_E(f^n(x)\in E,| R_E(x)=k)
=\rho_E(f^{n-k}(x)\in E)=e_{n-k}
$$
so that the $e_n$'s would satisfy the recurrence equation
$$
e_n-p_ne_0-p_{n-1}e_1-\cdots -p_1e_{n-1}=\cases{1, &for $n=0$,\cr
                                                0, &for $n>0$,\cr}
\eqno(A.11)
$$
where $p_n\equiv\rho_E(E_n)$. This would make $e_0,e_1,e_2,\dots$
the {\sl renewal sequence} 
associated to the sequence $p_1,p_2,\dots$.
It has been observed [Fe2] (see also [Ki]) that any renewal sequence, that is any sequence
generated as in (A.11) with $p_1,p_2\dots$ satisfying $p_n\geq 0$ and
$\sum p_n \leq 1$, can arise as the diagonal transition probabilites corresponding to a
given state in some Markov chain. In our case, 
a Markov chain which does the job is precisely that discussed
in Section 2, with the $p_n$'s as above and $e_n=p_{11}^n$.
Indeed, it is not difficult to realize that 
the Markov chain in question
is that with transition probabilities 
$$p_{ij}=\rho\, (A_i\cap f^{-1}A_j)/\, \rho(A_i)
\eqno(A.12)
$$ and stationary
distribution $\pi_i = \rho (A_i)$, where the sets $A_i$ are defined
in (A.4).

\noindent
We point out that under the supposition made 
above this Markov chain would be isomorphic (mod 0) 
to the iteration process $\{f^n(x)\}$.
On the other hand, in general the $R_E \circ f_E^{k}$ 
are not mutually independent 
and we are then led to call the above Markov chain the 
{\sl Markov approximation of
the dynamical system $(X,\rho,f)$ w.r.t. the reference set $E$}.
Leaving any further detail of this approximation procedure to 
be discussed elsewhere [Is1], in particular the question of the
choice of the reference set $E$ and that of the
``proximity" of $(X,\rho,f)$ and its Markov approximation 
(see [Che] where this and related questions for a closely
related approximation scheme have been
dealt with in a far reaching way), 
we are now going to discuss a simple example 
(modelling temporal intermittency)
where such an approximation is ``exact", in that
it is isomorphic 
to the dynamical system itself.
\noindent

\vsni
\noindent
{\bf Example.}
The Markov chain $P$ studied in Section 2 is isomorphic (mod $0$)
to the iteration process of the
piecewise affine `intermittent' map 
$f:[0,1]\to [0,1]$ given by
$$
f (x) =\cases{  (x-d_1)/ \alpha_1, &if $d_1\leq x\leq d_0$ \cr
        d_{i-1} + (x-d_i)/ \alpha_i, &if $d_i\leq x< d_{i-1},
\,\, i\geq 2$ \cr 
}\eqno(A.13)
$$
\noindent
Here the numbers $d_i=\sum_{l>i} p_l$ 
are supposed to be all distinct,
and $\alpha_i=p_i/p_{i-1}$, $i\geq 1$ (with $p_0=1$).
In what follows we shall always assume that $\sum d_i < \infty$.
The partition ${\cal A}$ of $[0,1]$
into the intervals $A_n=[d_{n},d_{n-1}]$, $n\geq 1$ is a Markov
partition for $f$.

\noindent
This map is named `intermittent' for,
if $\limsup \alpha_i =1$, then $f$ can be viewed as a
piecewise affine approximation of a piecewise smooth transformation
of $[0,1]$ which 
is expanding everywhere but at the 
fixed point in the origin, 
where the derivative is equal to one.

\noindent
 Let $\O$, $\pi$ be as in the example
of Section 2. 
One then sees that the map $\phi : \O \to [0,1]$ defined by:
$\phi (\o)=x$ according to $f^j(x)\in A_{\o_j}$, $j\geq 0$,
is a bijection between $\O$ and the residual set of points 
in $(0,1]$ which are not preimages of $1$ w.r.t. the map $f$.
Clearly $\phi$ conjugates $f$ with the shift $\tau$
on $\O$. Moreover, let $\mu$ be the $\tau$-invariant 
Markov probability measure on $\O$
defined in (1.8) (with $\pi$ and $P$ as above). 
Then $\rho =\mu \circ \phi^{-1}$ is $f$-invariant and
it is easy to see that 
the $p_{ij}$'s are as in $(A.4)-(A.12)$ with $E=[d_1,d_0]$.
Finally, if $f$ is the
piecewise affine approximation of a smooth transformation
of $[0,1]$ having a tangency at $x=0^+$ of order $1+1/\eta$,
with $\eta > 0$, then 
$p_i\sim i^{-(1+\eta)}$ and hence
$\alpha_i \sim  1 - (1+\eta)/ i$. Thus, in order to
have $\sum d_i < \infty$ it is necessary that $\eta > 1$, 
and the corresponding Markov chain $P$ has
ergodic degree $d =\eta - 1$.

\noindent
Let us consider the Perron-Frobenius operator
$M:L^1([0,1],dx)\to L^1([0,1],dx)$ which satisfies
$$
\int_0^1\, u\circ f^n(x) \, v(x) \, dx
=\int_0^1\, u(x) \, M^nv(x) \, dx\eqno(A.14)
$$
for all pairs $u,v\in L^1$. Note
that the space $\ell_1(S,p)$ of
vectors $u:\N\to \R$ such that
$$
\Vert u \Vert_{1,p}:=\sum_{i\in S}|u_i|\, p_i < \infty
$$
is left invariant by the operator $M$, 
which takes on the matrix representation
$$
M(i,j)={p_i\over p_j}\, P(i,j),\qquad i,j\geq 1.\eqno(A.15)
$$
The eigenequation $M\, h = h$ has a solution $h\in \ell_1(S,p)$
given by 
$h_i = h_1\, p_1\, d_{i-1}/ p_{i}$, and
the vector $p$
satisfies $M^*p=p$. Therefore,
recalling that $\pi_i = \pi_1 \, d_{i-1}$, and
putting $h_1= \pi_1/ p_1$, we get
$\pi_i = h_i\, p_i$.
One then sees that the vector $h \in \ell_1(S,p)$ corresponds to the 
(locally constant) density
of the absolutely continuous $f$-invariant probability measure 
$\rho (dx) = h(x)\, dx$, with $h\in L^1([0,1],dx)$ and $h(x)\equiv h_i$ for
$d_{i}\leq x < d_{i-1}$. Observe that $\rho (A_i)=\pi_i$.
Now, using (A.14) we find
$$
\rho (u\circ f^n \, v ) - \rho (u)\, \rho(v)=
\int_0^1\, u (x)\, 
[\, (M^n\, v h)(x) -  \rho (v)\, h(x)] \, dx\eqno(A.16)
$$
Suppose that $u$ and $v$ are bounded $L^1$-functions
taking constant values $u_i$ and $v_i$ on the elements $A_i$ of
the Markov partition ${\cal A}$. We shall denote 
by ${\bf u}=(u_i)_{i=1}^{\infty}$ and 
${\bf v}=(v_i)_{i=1}^{\infty}$ the 
corresponding vectors in $\ell_{\infty}(\N)$.
Using (A.15), (A.16) and the above observations
we get (the notation is as in the proof of 
Corollary 1),
$$
\rho (u\circ f^n \, v ) - \rho (u)\, \rho (v)
=({\overline {\pi {\bf v}} } \,
P^n-({\pi}\cdot {\bf v})\, \pi) \cdot {\bf u}.\eqno(A.17)
$$
Now set $u_{\infty}=\limsup u_i$, $v_{\infty}=\limsup v_i$,
and suppose that $u_{\infty}\neq 0$ or $v_{\infty}\neq 0$.
Then, setting
${\hat {\bf u}}={\bf u}-u_{\infty}{\bf 1}$ and 
${\hat {\bf v}}={\bf v}-v_{\infty}{\bf 1}$
have that
$(\pi \cdot {\hat {\bf u}}) (\pi \cdot {\hat {\bf v}})\neq 0$
provided $\pi \cdot {\bf u}\neq u_{\infty}$ and
$\pi \cdot {\bf v}\neq v_{\infty}$. Moreover 
$\limsup {\hat u}_i=0$ and $\limsup {\hat v}_i=0$.
On the other hand we plainly have 
$({\overline {\pi {\hat {\bf v}}} } \,
P^n-({\pi}\cdot {\hat {\bf v}})\, \pi) \cdot {\hat {\bf u}}=
({\overline {\pi {\bf v}} } \,
P^n-({\pi}\cdot {\bf v})\, \pi) \cdot {\bf u}$.
We then see that the conditions
$\rho (u)\equiv \pi \cdot {\bf u}\neq u_{\infty}$ and
$\rho (v)\equiv \pi \cdot {\bf v}\neq v_{\infty}$ 
are equivalent to the conditions
$u_i=o (1)$ and $\nu_i=o (\pi_i)$ 
(along with $(\pi \cdot {\bf u}) (\nu \cdot 1)\neq 0$) 
assumed in the last statement of Theorem 2,
with
the identification 
$\nu = {\overline {{\pi} {\bf v}} } /{\pi}\cdot {\bf v}$.

\noindent
The following result is now a direct 
consequence of Theorem 2
(for related results see [Is2], [LSV], [Mo]; see also [Yo], [Is1] and [Sa]
for more general approaches dealing with smooth maps): 

\vsni
\noindent
{\bf Corollary 2.} {\sl Let $f:[0,1]\to [0,1]$ be as in
(A.13) and assume that $\alpha_i \sim  1 - {(1+\eta)/ i}$
for some $\eta >1$.
Then, for any pair of bounded $L^1$-functions $u,v:[0,1]\to \R$,
locally constant on the Markov partition ${\cal A}$, there is a 
positive constant $C=C(u,v)$ such that, for $n$ large enough,
$$
|\,\rho (u\circ f^n \, v ) - \rho (u)\, \rho (v)\, | \leq
C\, n^{-(\eta-1)}.
$$
Assume furthermore that $\rho (u)\neq u_{\infty}$ and
$\rho (v)\neq v_{\infty}$. Then we have}
$$
\rho (u\circ f^n \, v ) - \rho (u)\, \rho (v) \sim 
C\, n^{-(\eta-1)}.
$$
We conclude with a final remark. 
From the proof of Theorem 2 it follows that
if the conditions $\rho (u)\neq u_{\infty}$ and
$\rho (v)\neq v_{\infty}$ 
are violated, then cancellations may take place 
to accelerate the convergence rate. 
A trivial example is obtained by taking $u,v$ constant on $[0,1]$. 
Conversely, one may argue as follows: 
take $a>0$ and let $t_a(x)$ be the first entrance time into 
the set $[a,1]$. 
When an orbit falls in a small (compared to $a$) neighbourhood
of $0$ it stays there for a time which can be 
arbitrarily large before reaching again $[a,1]$. 
More precisely, from the above discussion one readily finds that,
under the assumptions of Corollary 2,
$$
\rho \{x\in [0,1] \, :\, t_a(x)>n\}\sim C(a)  \, 
n^{-(\eta-1)}.
$$
Thus, if the condition is satisfied, namely if
the average value of the test functions is reached
away from the origin, then
the term $\rho (u\circ f^n \, v )$ cannot approach its asymptotic
value $\rho (u)\, \rho (v)$ at a rate faster 
than that given by the statistics
of first entrance times given above. 

\vskip .5cm
\noindent
{\sl Acknowledgements}: I would like to thank Lai-Sang Young
for interesting conversations at the origin of this research and
the referee for a constructive criticism and 
several valuable suggestions.

\vskip 0.5cm
\noindent
{\bf References.}
\vsni
\noindent

\item{[A]} M Amici: Honour thesis: Propriet\`a statistiche di 
alcuni processi di rinnovamento,
{\sl the University of Camerino}, 2001.

\item{[Che]} N Chernov: Limit theorems and markov approximations for
chaotic dynamical systems, {\sl Probab. Theory Relat. Fields} {\bf 101},
(1995) 321-362.

\item{[Chu]} K L Chung: {\it Markov chains with stationary
transition probabilities},
Springer-Verlag Berlin Heidelberg New York 1967.

\item{[De]} C Derman: 
Some contributions to the theory of denumerable Markov chains,
{\sl TAMS} {\bf 73} (1955), 471-486.  

\item{[Fe1]} W Feller: 
Fluctuation theory of recurrent events,
{\sl TAMS} {\bf 67} (1949), 99-119.  

\item{[Fe2]} W Feller: 
{\sl An Introduction to Probability Theory and Its
Applications}, Volume 2,
J.Wiley and Sons, New York 1970.

\item{[FMM]} G Fayolle, V A Malyshev and M V Menshikov: 
{\it Topics in the constructive theory of countable markov
chains},
Cambridge University Press, Cambridge 1992. 
 
\item{[H1]} T E Harris: 
First passage and recurrence distribution,
{\sl TAMS} {\bf 73} (1952), 471-486.  

\item{[H2]}
G H Hardy: {\it Divergent series},  Oxford at the Calrendon Press 1949.

\item{[Is1]} S Isola: On systems with finite ergodic degree, 
Preprint 2001.

\item{[Is2]} S Isola: Renewal sequences and intermittency, 
{\sl J. Stat. Phys.} {\bf 97} (1999), 263-280.

\item{[Ki]} J F C Kingman: {\it Regenerative phenomena}, John Wiley,
1972.

\item{[KSK]} J G Kemeny, J Snell and A W Knapp: 
{\it Denumerable Markov Chains}, Van Nostrand, Princeton, 1966.

\item{ [LSV]}  A Lambert, S Siboni and S Vaienti: {Statistical
properties of a non-uniformly hyperbolic map of the interval},
{\sl J. Stat. Phys.} {\bf 72} (1993), 1305-1330.

\item{[Mo]} M Mori: {On the intermittency of a piecewise
linear map}, {\sl Tokyo J. Math.} {\bf 16} (1993), 411-428.  

\item{[Pi]} J W Pitman: Uniform rates of 
convergence for Markov chains transition probabilities, 
{\sl Z. Wahrscheinlichkeitstheorie verw. Gebeite},
{\bf 29} (1974), 193-227. 
 
\item{[Pop]} N N Popov: On the rate of convergence for countable 
Markov chains, {\sl Theory Prob. Appl.} {\bf 24} (1978), 401-405.

\item{[Pos]}
A G Postnikov: {\sl Tauberian Theory and its Applications}, 
Proceedings of the Steklov Institute
of Mathematics, 1980, Issue 2. 

\item{ [Ro]}  B A Rogozin: An estimate of the remainder term
in limit theorems of renewal theory, 
{\sl Theory Prob. Appl.} {\bf 18} (1973), 662-677. 

\item{ [Ru]} D Ruelle: {\sl Thermodynamic Formalism}, 
 Addison-Wesley Publ. Co. 1978.

\item{ [Sa]}  O Sarig: Subexponential decay of correlations, 
Preprint 2001. 

\item{ [Se]}  B A Sevast'yanov: Renewal theory, 
{\sl J. Soviet Math.} {\bf 4} (1975), n.3. 

\item{[TT]} P Tuominen and R L Tweedie: Subgeometric 
rates of convergence of $f$-ergodic markov chains, 
{\sl Adv. Appl. Prob.} {\bf 26} (1994),
775-798. 

\item{[VJ]} D Vere-Jones: On the spectra of some linear operators
associated with queueing systems, {\sl Z. Wahrsch.} {\bf 2} (1963),
12-21. 

\item{ [Wa]} X J Wang:
{Statistical physics of temporal intermittency}, 
{\sl Phys. Rev.} {\bf A40} (1989), 6647.

\item{[Yo]} L S Young: Recurrence times and rate of mixing,
 {\sl Isr. J. Math.} {\bf 110} (1999),  153-188.

\end